\documentclass[12pt]{article}
\usepackage{amsmath,amssymb,theorem}
\usepackage{graphicx}
\usepackage{subfigure}
\usepackage{psfrag}
\usepackage{color}
\usepackage{enumerate}
\usepackage{epstopdf}
\usepackage[ruled, lined, linesnumbered]{algorithm2e}
\newtheorem{thm}{Theorem}[section]
\newtheorem{conj}[thm]{Conjecture}
\newtheorem{cor}[thm]{Corollary}
\theorembodyfont{\rmfamily}
\def\pf{\bigskip\noindent {\bf Proof.}~~}

\def\dfn#1{{\sl #1}}

\def\less{\backslash}

\def\pf{\bigskip\noindent {\bf{Proof.}}~~}

\def\mytextindent#1{\indent\llap{#1\enspace}\ignorespaces}

\def\myitem{\par\hangindent\parindent\mytextindent}

\title{Antimagic orientations of even regular graphs}
\author{Tong Li$^1$, Zi-Xia Song$^2$ , Guanghui Wang$^1$\thanks{Corresponding author.  \newline Email addresses:  tongli121@163.com (T. Li), Zixia.Song@ucf.edu (Z-X. Song),  ghwang@sdu.edu.cn (G. Wang), dlyang120@163.com (D. Yang), cqzhang@math.wvu.edu (C-Q. Zhang)},  Donglei Yang$^1$, Cun-Quan Zhang$^3$ \\
{ $^1$ Department of Mathematics}\\
{ Shandong  University, Jinan,  China}\\
{ $^2$ Department  of Mathematics}\\
{University of Central Florida}\\
{ Orlando, FL32816, USA}\\
{ $^3$ Department of Mathematics}\\
{ West Virginia University}\\
{ Morgantown, WV26506, USA}\\
}

\begin{document}
\maketitle
\begin{abstract}
A \dfn{labeling} of a digraph $D$ with $m$ arcs is a bijection from the set of arcs of $D$ to $\{1, \ldots, m\}$. A labeling of $D$ is  \dfn{antimagic} if no two vertices in $D$ have the same  vertex-sum,  where the  vertex-sum  of a vertex $u\in V(D)$  for a labeling is the sum of labels of all arcs entering $u$ minus the sum of labels of all arcs leaving $u$.
Motivated by the conjecture of Hartsfield and Ringel from  1990 on antimagic labelings of graphs, Hefetz, M\"{u}tze, and Schwartz [On antimagic directed graphs, J Graph Theory 64 (2010) 219--232] initiated the study of antimagic labelings of digraphs, and  conjectured that every connected graph  admits an antimagic orientation, where  an orientation $D$ of a graph $G$ is \dfn{antimagic} if  $D$ has an antimagic labeling. It remained unknown whether  every disjoint union of cycles admits an antimagic orientation.
 In this paper, we first  answer this question in the positive  by proving that every  $2$-regular graph has an antimagic orientation. We then show that  for any integer $d\ge2$, every connected,  $2d$-regular graph has an antimagic orientation.   Our technique is new.

\end{abstract}
\bigskip
\noindent {\textbf{Keywords}: regular graph,  antimagic labeling, antimagic orientation

\baselineskip 18pt
\section{Introduction}
All graphs   in this paper are finite and simple.   For a graph $G$, we use $|G|$ and $e(G)$ to denote the number of vertices and edges of $G$, respectively.   An \dfn{antimagic labeling} of a graph $G$ is a bijection from $E(G)$ to $\{1,2,..., e(G)\}$ such that for any distinct vertices $u$ and $v$, the sum of labels on edges incident to $u$ differs from that for edges incident to $v$. A graph $G$ is \dfn{antimagic} if it has an antimagic labeling. Hartsfield and Ringel~\cite{NG} introduced antimagic labelings in 1990 and conjectured that every connected graph  other than $K_2$  is  antimagic.  The most significant progress on this problem is a result of Alon, Kaplan, Lev, Roditty, and Yuster~\cite{Alon}, which states that there exists an absolute   constant  $c$ such that  every   graph on $n$ vertices with minimum degree at least $c\log n$  is antimagic.  Eccles~\cite{E} recently improved this result by showing that there exists an absolute constant $c_0$ such  that if   $G$  is a graph with average
degree at least $c_0$, and $G$ contains no isolated edge and at most one isolated vertex, then $G$ is
antimagic.    Cranston~\cite{DC} proved that any $d$-regular bipartite graph with   $d\ge2$ is antimagic. For nonbipartite regular graphs, Cranston, Liang,  and Zhu~\cite{DYZ} proved that every odd regular graph is  antimagic, and later Chang, Liang, Pan, and Zhu~\cite{CLPZ} proved that every even regular graph is  antimagic.  For more information on antimagic labelings of graphs and related labeling problems, see the  recent  informative survey~\cite{JAG}.\medskip


Motivated by antimagic labelings of graphs, Hefetz, M\"{u}tze, and Schwartz~\cite{DTJ} initiated the study of antimagic labelings of digraphs. For a positive integer $k$, we define $[k]:=\{1,2,\dots, k\}$. Let  $D$ be a digraph. We use $A(D)$ and $V(D)$ to denote the set of arcs and vertices of $D$, respectively.  A \dfn{labeling} of  $D$ with $m$ arcs is a bijection from $A(D)$ to $[m]$. A labeling of $D$ is  \dfn{antimagic} if no two vertices in $D$ have the same  vertex-sum,  where the  vertex-sum  of a vertex $u\in V(D)$  for a labeling is the sum of labels of all arcs entering $u$ minus the sum of labels of all arcs leaving $u$.  A digraph $D$ is \dfn{antimagic} if it has an antimagic labeling. A graph $G$ has an \dfn{antimagic orientation} if an orientation  of  $G$ is  antimagic.  Hefetz, M\"{u}tze, and Schwartz~\cite{DTJ} raised the questions
`` Is every orientation of any connected  graph antimagic?" and ``Does every  graph admit an antimagic orientation?".  Except for $K_{1,2}$ and $K_3$, no other counterexamples to the first question are known.  They proved  an analogous result of Alon, Kaplan, Lev, Roditty, and Yuster~\cite{Alon} that there exists an absolute constant $c$ such that every orientation of any  graph on $n$ vertices with minimum degree at least $c\log n$ is antimagic. They also showed that every orientation of   star $S_n$ with  $n\ne 2$ is antimagic; every orientation of  wheel $W_n$ is antimagic; and every orientation of $K_n$ with  $n\ne3$ is antimagic.  For the second question, they prove the following. \medskip

\begin{thm}[\cite{DTJ}]\label{odd} For any integer $d\ge1$,
\myitem{(a)} every   $(2d-1)$-regular graph admits an antimagic orientation.
\myitem{(b)} every connected, $2d$-regular graph $G$ admits an antimagic orientation if  $G$ has a matching covers all but at most one vertex of $G$.
\end{thm}

Hefetz, M\"{u}tze, and Schwartz~\cite{DTJ} asked  whether it is true that every orientation of any connected  graph  on at least four vertices is antimagic. They also pointed out that
 ``It seems hard to discard any of the two conditions in Theorem~\ref{odd}(b), that is connectedness and having a matching that covers all vertices but at most one. In fact, we do not even know if every disjoint union of cycles admits an antimagic orientation."  They proposed the following conjecture.

\begin{conj}[\cite{DTJ}]\label{conj1}
Every connected graph admits an antimagic orientation.
\end{conj}

 Recently, Shan and Yu~\cite{SY} proved that  Conjecture~\ref{conj1} holds  for  biregular bipartite graphs. It remained unknown whether  every $2$-regular graph, that is, every disjoint union of cycles, has an antimagic orientation.
  In this paper, we first answer this question in the positive by proving that  every  $2$-regular graph admits an antimagic orientation. We  then  prove that for any integer $d\ge2$, every  connected, $2d$-regular graph admits an antimagic orientation.
It turns out that finding an antimagic orientation of a $2$-regular graph is,   indeed,   a bit more complicated  than finding an antimagic orientation of an odd regular graph (see Theorem 1.3  in \cite{DTJ}) or a connected,  even regular graph (see Theorem~\ref{th2} below). Our  technique is new and proofs of both results are neat. \medskip

We need to introduce  more notation.  A closed walk in a graph  is an \dfn{Euler tour} if
it traverses every edge of the graph exactly once.  
The following is a result of Euler.

\begin{thm}[Euler 1736]\label{euler}
A connected graph admits an  Euler tour if and only if every vertex has even degree.
\end{thm}

Let  $D$ be an orientation of a graph $G$ with $m$ edges. For any   labeling $c: A(D)\rightarrow [m]$ of $D$ and any vertex $u\in V(D)$,  we use $s_D(u)$,  or simply $s(u)$ when there is no confusion,  to denote the vertex sum of $u$ for the labeling  $c$.

\section{Antimagic orientations of $2$-regular graphs}

In this section we study  antimagic orientations of $2$-regular graphs.  It remained unknown whether  every $2$-regular graph has an antimagic orientation.   We answer this question in the positive below.

\begin{thm}\label{th1}
Every  $2$-regular graph admits an antimagic orientation.
\end{thm}

\pf  Let $G$ be a $2$-regular graph on $n$ vertices. Then $e(G)=n$ and every component of $G$ is a cycle. Let $C_1$, $\dots$, $C_s$, $C_{s+1}$, $\dots$, $C_{s+t}$ be all distinct components of $G$ such that $C_1, \dots, C_s$ are odd cycles  and $C_{s+1}, \dots, C_{s+t}$ are even cycles, where $|C_1|\le \cdots\le |C_s|$ and   $|C_{s+1}|\le\cdots\le |C_{s+t}|$. For any $i\in[s+t]$, we may assume that $C_i$ has vertices $v_{i,1}, v_{i,2}, \dots, v_{i,r_i}$ in order, where $r_i:=|C_i|$.  We first find an orientation $D$ of $G$. Let $M_1$ and $M_2$
be two disjoint perfect matchings of    $G\less \{v_{1,1}v_{1,r_1}, \dots, v_{s,1}v_{s, r_s}\}$ such that $v_{1,1}v_{1, 2}, \dots, v_{s+t, 1}v_{s+t, 2}\in M_1$. Let $D$ be the orientation of $G$ by directing the  edges  of $G$ as follows: for all $i\in[s]$, orient every edge $v_{i, 1}v_{i, r_i}$ in $\{v_{1,1}v_{1, r_1}, \dots, v_{s, 1}v_{s, r_s}\}$ from   $v_{i, r_i}$ to $v_{i, 1}$;  then  for all $i\in[s+t]$,  orient every edge $v_{i,j}v_{i, j+1}\in M_1$  from    $v_{i, j}$ to $v_{i, j+1}$;  and  every edge $v_{i, j}v_{i, j+1}\in M_2$ from $v_{i, j+1}$ to  $v_{i, j}$, where  all arithmetic on the index $j+1$ in $v_{i, j+1}$ for each cycle $C_i$  here and henceforth is done modulo $r_i$. Clearly, $D$ is an orientation of $G$. Let $D_o$ be the above orientation of the odd cycles $C_1, \dots, C_s$. \medskip

We next find  a labeling $c: A(D)\rightarrow [n]$  of $D$ such that $c$, together with $D$,  is a desired antimagic orientation of $G$.   Let $n_0: =0$ and $n_i:=r_1+\cdots+r_i$ for all $i\in [s+t]$.  We first find a bijection $c_e: A(C_{s+1})\cup\cdots\cup A(C_{s+t})\rightarrow \{n_s+1, \dots, n\}$.  For any $i\in\{s+1, \dots, s+t\}$, let $c_e(v_{i, j}v_{i, j+1})=n_{i-1}+j$ for all $j\in [r_i-2]$, $c_e(v_{i, r_{i-1}}v_{i, r_i})=n_{i-1}+r_i$, and $c_e(v_{i, r_i}v_{i,1})=n_{i-1}+r_i-1$.
We next find a labeling   $c_o: A(D_o)\rightarrow [n_s]$ of $D_o$ with  $s(v_{i, 1})=-i$  for all $i\in [p]$,  and $s(v_{p+j, 1})=j$ for all  $j\in \{1, \dots, \lceil\frac s2\rceil\}$, where $p= \lfloor\frac s2\rfloor$.  Let $c_o(v_{1, r_1}v_{1,1})=1$ and $c_o(v_{{p+1}, r_{p+1}}v_{{p+1}, 1})=n_s$.  Then  $c_o(v_{1,1}v_{1,2})=c_o(v_{1, r_1}v_{1,1})-s(v_{1, 1})=2$, and $c_o(v_{{p+1}, 1}v_{{p+1}, 2})=c_o(v_{{p+1}, r_{p+1}}v_{{p+1}, 1})-s(v_{p+1, 1})=n_s-1$. We then label the remaining edges of $C_1, \dots, C_s$ recursively as
depicted in Algorithm~\ref{label1} on the next page, where the edges of $C_1, \dots, C_p$ are labelled from line 1 through line 9, and the edges of $C_{p+1}, \dots, C_s$ are labelled from line 10 to line 18.  Let $c$ be obtained from $c_o$ and $c_e$, that is, label the arcs in $D$ as they are labeled under $c_o$ and $c_e$. Clearly, $c$ is a labeling of $D$. \medskip

\begin{algorithm}[htb]\label{label1}
\SetAlgoLined
\KwData{Odd cycles $C_1, \dots, C_s$ with the given orientation $D_o$,    $s(v_{i, 1})=-i$ for all $i\in[p]$,  $s(v_{j, 1})=j-p$ for all $j\in\{p+1, \dots, s\}$,  $c_o(v_{1, r_1}v_{1,1})=1$, $c_o(v_{1,1}v_{1,2})=2$,  $c_o(v_{{p+1}, r_{p+1}}v_{{p+1}, 1})=n_s$, and $c_o(v_{{p+1}, 1}v_{{p+1}, 2})=n_s-1$}
\KwResult{An antimagic labeling of $D_o$}

\For{$i=2$ to $p$ }{
Assign the smallest unused number, say $\alpha$,  in $[n_s]$ to the edge entering   $v_{i,1}$\;

Assign the value $\alpha-s(v_{i,1})$ to the edge leaving the vertex $v_{i,1}$\;
}

Set  $A$ to be the set of edges in $G$  incident with $v_{i,1}$ for all   $i\in [p]$, and  set $A^*$ to be $A$\;

\While{$A\neq E(C_1)\cup\cdots\cup E(C_p)$}{

Assign the smallest unused number in $[n_s]$ to the edge $e\in (E(C_1)\cup\cdots\cup E(C_p))\less A$  which is adjacent to the edge $e^*\in A^*$ with $c_o(e^*)$ the smallest among the labels on the edges in $A^*$\;

Set $A$ to be $A\cup\{e\}$, and set $A^*$ to be $(A^*\less e^*)\cup\{e\}$\;
}

\For{$j=p+2$ to $s$ }{
Assign the largest unused number, say $\beta$,  in $[n_s]$ to the edge entering  $v_{j,1}$\;

Assign the value $\beta-s(v_{j,1})$ to the edge leaving $v_{j,1}$\;
}

  Set  $B$ to be the set of edges in $G$ incident with $v_{j,1}$ for all $j\in \{p+1, \dots, s\}$,  and set $B^*$ to be $B$\;
\While{$B\neq E(C_{p+1})\cup\cdots\cup E(C_s)$}{

Assign the largest unused number in $[n_s]$ to the edge $e\in (E(C_{p+1})\cup\cdots\cup E(C_s))\less B$ which is adjacent to the edge $e^*\in B^*$ with $c_o(e^*)$ the largest among the labels on the edges in $B^*$  \;

Set $B$ to be $B\cup\{e\}$, and set $B^*$ to be $(B^*\less e^*)\cup\{e\}$\;

}

\caption{Label the edges of $C_{1}, \dots, C_s$}
\end{algorithm}

 It remains to verify that  $c$ is antimagic.
For any $i\in\{s+1, \dots, s+t\}$, we see that  $s(v_{i,1}) =-(2n_{i-1}+r_i)$, $s(v_{i,j}) =(-1)^{\delta_j}(2n_{i-1}+2j-1)$ for all $j\in\{2, 3, \dots, r_i-2, r_i\}$, and $s(v_{i,r_i-1}) =-(2n_{i-1}+2r_i-2)$, where $\delta_j=0$ if $j$ is even and  $\delta_j=1$ if $j$ is odd. Clearly,  no two vertices of $C_{s+1},\dots, C_{s+t}$ have the same vertex-sum under $c$.  Thus $c$ is an antimagic labeling of $D$ if $s=0$. So we may assume that $s\ge1$.  Next, for any $u\in V(C_1)\cup\dots\cup V(C_s)$ and $v\in V(C_{s+1})\cup\dots\cup V(C_{s+t})$,  we see that $|s(u)|\le 2n_s-1$ and $|s(v)|\ge  2n_s+3$. Thus $s(u)\ne s(v)$. To show that $c$ is antimagic, it suffices to show that $c_0$ given in  Algorithm~\ref{label1} is an antimagic labeling of $D_o$. We do that next. \medskip

Let $X:=\{v_{1,1}, v_{2,1}, \dots, v_{s,1}\}$.  By the choice of $c_o$, no two vertices in $X$ have the same vertex-sum under $c_o$. Furthermore, for any  $u\in (V(C_1)\cup\dots\cup V(C_s))\less X$, by the orientation $D$ of $G$, $|s(u)|=a+b$ for some distinct integers  $a, b\in [n_s]$.  According to line 1 through line 9 in Algorithm~\ref{label1}, $\alpha\le n_p$ and for any  $u\in (V(C_1)\cup\dots\cup V(C_p))\less X$, either  $-(2n_p-1)\le s(u)\le -(p+2)$ or $p+3\le s(u)\le 2n_p-1$.   Similarly,  according to line 10 through line 17 in Algorithm~\ref{label1}, $\beta\ge n_p+1$ and   for any  $v\in (V(C_{p+1})\cup\dots\cup V(C_s))\less X$, either  $ -(2n_s-1)\le s(v)\le -(2n_p+3)$ or $2n_p+3\le s(v)\le 2n_s-1$.  It follows that no vertex in $X$ has the same vertex-sum as any vertex in   $ (V(C_1)\cup\dots\cup V(C_s))\less X$; and no vertex in $ (V(C_1)\cup\dots\cup V(C_p))\less X$ has the same vertex-sum as any vertex in  $(V(C_{p+1})\cup\dots\cup V(C_s))\less X$. It remains to show that
 no two vertices in $ (V(C_1)\cup\dots\cup V(C_p))\less X$ (resp.  $(V(C_{p+1})\cup\dots\cup V(C_s))\less X$ ) have the same vertex-sum. We only verify below  that  no two vertices in $ (V(C_1)\cup\dots\cup V(C_p))\less X$  have the same vertex-sum, because by the choice of labels stated in the line 7 and line 16 in Algorithm~\ref{label1}, a similar argument can be applied to prove that no two vertices in $(V(C_{p+1})\cup\dots\cup V(C_s))\less X$ have the same vertex-sum.   \medskip

By Algorithm~\ref{label1}, we observe the following.

\begin{enumerate}[(a)]
\item   For any $i, k\in[p]$ with $i<k$,  $c(v_{k, r_k-j+1}v_{k, r_k-j})<c(v_{i, r_i-j}v_{i, r_i-j-1})< c(v_{k, r_k-j}v_{k, r_k-j-1})$ for any  $j=0, 1, \dots, \lfloor\frac{r_i}2\rfloor-2$,  and  $c(v_{i, j}v_{i, j+1})<c(v_{k, j}v_{k, j+1})$ for any  $j= 1, \dots, \lfloor\frac{r_i}2\rfloor$. 
\item  For any  $i\in[p]$,  $c(v_{i, r_i-j}v_{i, r_i-j-1}) > c(v_{i, j+1}v_{i, j+2})$ for any $j=0 \dots, \lfloor\frac{r_i}2\rfloor-2$.
\end{enumerate}

Suppose for a contradiction that there exist two distinct vertices  $u,w$ in $ (V(C_1)\cup\dots\cup V(C_p))\less X$  such that $s(u)=s(v)$. We may assume that $u\in V(C_i)$ and $w\in V(C_k)$ for some $i, k\in[p]$. Clearly, $i\ne k$. 
We may further assume that $i<k$,   $s(u)=a+b$,  and $s(w)=x+y$ for some distinct integers $a,b, x,y\in [n_p]$ with $a<b$ and $x<y$. If  $a<x$ and $u\neq v_{i, \lfloor\frac{r_{i}}{2}\rfloor}$, then by line 7 in Algorithm~\ref{label1} , $b<y$, contrary to the fact that $a+b=x+y$. If  $a>x$ and $u\neq v_{i, \lfloor\frac{r_{i}}{2}\rfloor}$, then by (a) and line 7 in Algorithm~\ref{label1} , $b>y$, a contradiction. If  $a<x$ and $u=v_{i, \lfloor\frac{r_{i}}{2}\rfloor}$, since $c(v_{i, \lfloor\frac{r_{i}}{2}\rfloor+2}v_{i, \lfloor\frac{r_{i}}{2}\rfloor+3})<a$, then $c(v_{i, \lfloor\frac{r_{i}}{2}\rfloor+2}v_{i, \lfloor\frac{r_{i}}{2}\rfloor+3})<x$, by line 7 in Algorithm~\ref{label1} , $b<y$, also a contradiction. Thus $a>x$ and $u=v_{i, \lfloor\frac{r_{i}}{2}\rfloor}$, then $x<a<b<y$. Let $q:=\lfloor\frac{r_i}2\rfloor$ and $q^*:=\lfloor\frac{r_k}2\rfloor$.  Then $q\le q^*$ because $r_i\le r_k$. According to  Algorithm~\ref{label1}, we see that $u=v_{i, q+1}$, and $c(v_{i, r_i-q+2}v_{i, r_i-q+1})<x$. 
Suppose first that  $w=v_{k, j}$ for some $j\in [q^*]$.  By the orientation of  $G$ and the fact that $a>x$, we see that $j\le q-1$.   By (a), $c(v_{i, r_i-q+2}v_{i, r_i-q+1})>c(v_{k, r_k-j+2}v_{k, r_k-j+1})$  because $j\le q-1$.
By  (b), $c(v_{k, r_k-j+2}v_{k, r_k-j+1})>x$.
 It follows that  $c(v_{i, r_i-q+2}v_{i, r_i-q+1})>x$, contrary to the fact that $c(v_{i, r_i-q+2}v_{i, r_i-q+1})<x$.  Thus $w=v_{k, r_k-j}$ for some $j\in [q^*]$.  Since $a>x$, by the orientation of  $G$ and (a), we see that $j\ge q$. By (b), $c(v_{i, r_i-q+2}v_{i, r_i-q+1})>c(v_{i, q-1}v_{i, q})$.  Since $a>x$,  by line 7 in Algorithm~\ref{label1}, $c(v_{i, q-1}v_{i, q})> c(v_{k, r_k-j+2}v_{k, r_k-j+1})$.
 It follows that  $c(v_{i, r_i-q+2}v_{i, r_i-q+1})>c(v_{k, r_k-j+2}v_{k, r_k-j+1})$, which is impossible because $j\ge p$.   
 \medskip

This completes the proof of Theorem~\ref{th1}.   \hfill\vrule height3pt width6pt depth2pt\\
\section{Antimagic orientations of  even regular graphs}

In this section we  first prove a result on   antimagic orientations of connected, $2d$-regular graphs, where $d\ge2$.

\begin{thm}\label{th2}
For any integer $d\ge2$, every  connected, $2d$-regular graph admits an antimagic orientation.
\end{thm}

\pf  For any integer $d\ge2$, let $G$ be a connected, $2d$-regular graph on $n$ vertices. 
By Theorem~\ref{euler},  let $C^*$ be an   Euler tour of $G$. We can regard $C^*$ as a cycle $C$ with $d\ge2$ copies of each vertex of $G$  on $C$. For each vertex $v$ in $G$, arbitrarily pick one of the $d$ copies of $v$ on $C$ as a real vertex and the remaining  $d-1$ copies of $v$ as imaginary vertices. Then $C$ has $n$ real vertices and $(d-1)n$ imaginary vertices. Let $V_R=\{v_1, v_2, \ldots,v_{n}\}$  and $V_I=\{u_1,u_2,\ldots,u_{(d-1)n}\}$ be the set of real vertices and imaginary vertices of $C$, respectively. Then $V(C)=V_R\cup V_I$. By renaming the vertices in $V_R$ if necessary, we label the vertices of $V_R$ on $C$ with $v_1, v_2, v_4, \dots, v_n, v_{n-1}, v_{n-3}, \dots, v_3$ in order when $n$ is even;  and $v_1$, $v_2$, $v_4$, $\dots$,  $v_{n-1}$, $v_n$, $v_{n-2}$, $\dots$, $v_3$ in order when $n$ is odd, as depicted in Figure~\ref{C}.  Two vertices $v_i, v_j\in V_R$  are \dfn{good pair on $C$}  if there exists a $(v_i, v_j)$-path $P_{i, j}$ along $C$ so that  either $v_iv_j\in E(C)$  or all the internal vertices of $P_{i, j}$  are imaginary vertices. Notice that such a path $P_{i,j}$ is unique for any good pair  $v_i, v_j\in V_R$.  We next find  an orientation $D$ of $C$. \medskip

\begin{figure}[htb]
\centering
\includegraphics[scale=0.4]{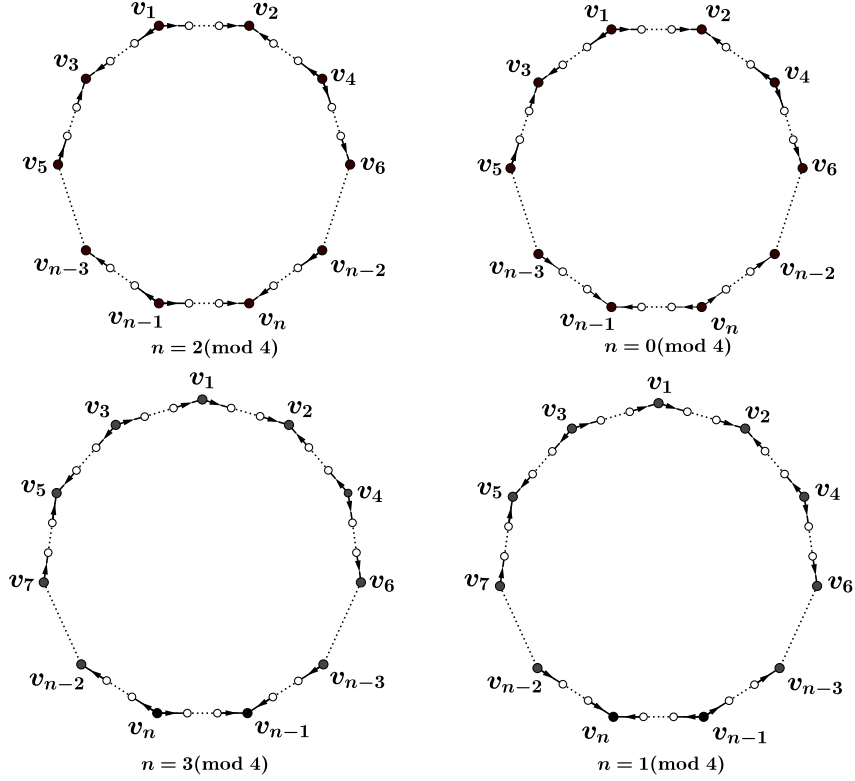}
\caption{Orientations of $C$ according to the parity of $n$, where real vertices and imaginary vertices of $C$ are denoted by  $\bullet$ and  $\circ$, respectively. }
\label{C}
\end{figure}

When $n$ is even,  set  $d_D^+(v_i)\in\{0,2\}$ for any $i\in [n]$, and $d_D^+(u_j)=1$ for any $j\in[(d-1)n]$ by  first  directing the  path $P_{1,2}$   from $v_1$ to $v_2$,   and   then paths $P_{2,4}$  from $v_4$ to $v_2$,  $\dots$,  and  finally $P_{5,3}$  from $v_5$ to $v_3$, and  $P_{3, 1}$ from $v_1$ to $v_3$.
When $n$ is odd,  set   $d_D^+(v_1)=1$,  $d_D^+(v_i)\in\{0,2\}$ for any $i\in \{2,3, \dots, n\}$, and $d_D^+(u_j)=1$ for any $j\in[(d-1)n]$ by  first  directing the  paths $P_{1,2}$   from $v_1$ to $v_2$  and  $P_{3,1}$ from $v_3$ to $v_1$, and then  paths $P_{2,4}$  from $v_4$ to $v_2$,  $\dots$,  and  finally $P_{5,3}$  from $v_3$ to $v_5$.  Orientations of $C$ for both cases are depicted in Figure~\ref{C}.  \medskip

\medskip

We need to  find a labeling $c: A(D)\rightarrow [dn]$ such that $c$, together with $D$,  is a desired antimagic orientation of $G$. Let $\ell_0$ be the length of $P_{1,2}$. For any $i\in [n-2]$, let $\ell_{i}$ be the length of $P_{i, i+2}$. Finally, let $\ell_{n-1}$ be the length of $P_{n-1, n}$.  Clearly,  $\sum_{j=0}^{n-1}\ell_j=dn$.  We define a bijection $c: A(D)\rightarrow [dn]$   as
stated in Algorithm~\ref{label2} below. \medskip

\begin{algorithm}[H]\label{label2}
\SetAlgoLined
\KwData{Cycle $C$ with the given orientation $D$}
\KwResult{A bijection  $c: A(D)\rightarrow [dn]$}

 Assign the numbers in $[\ell_0]$ to the edges of $P_{1,2}$ in the increasing order along the orientation of  $P_{1,2}$ \;
Set  $V=\{v_1, v_2\}$ \;
 \While{$V\neq \{v_1, \dots, v_{n}\}$}{

\For{$i=1$ to $n-2$ }{
 Assign the  numbers in $\{\ell_0+\cdots+\ell_{i-1}+1, \dots, \ell_0+\cdots+\ell_{i-1}+\ell_i\}$  to the edges of $P_{i, i+2}$ in the  increasing order along  the orientation of  $P_{i, i+2}$ \;

Set  $V$ to be $V\cup \{v_{i+2}\}$ \;
}{

 Assign the  numbers in $\{\ell_0+\cdots+\ell_{n-2}+1, \dots, \ell_0+\cdots+\ell_{n-2}+\ell_{n-1}\}$ to the edges of $P_{n-1, n}$  in the  increasing order along  the orientation of  $P_{n-1, n}$ \;
}
}

\caption{Label the arcs of $D$}
\end{algorithm}
\bigskip

By Algorithm~\ref{label2}, $s_D(u_j)=-1$ for all $j\in [(d-1)n]$.  Let $D^*$ be the corresponding orientation of $C^*$, that is, oriented each edge on $C^*$ as it is oriented on $C$. Clearly,  $D^*$ is an orientation of $G$.   It remains to verify that the bijection  $c$ given in Algorithm~\ref{label2} is an antimagic labeling of $D^*$.  We may assume that $V(G)=V_R$.
 For each $v_i\in V(G)$,  $s_{D^*}(v_i)=s_D(v_i)+(d-1)s_D(u_i^*)=s_D(v_i)-(d-1) $, where $u_i^*$ is one of the $d-1$ imaginary vertices of $v_i$. It suffices to show that  for any $v_i, v_j\in V_R$ with $i\ne j$, $s_D(v_i)\neq s_D(v_j)$. According to Algorithm~\ref{label2}, when $n$ is even, we see that $s_D(v_1)=-\ell_0-2$;   $s_D(v_n)=-2(\ell_0+\cdots+\ell_{n-3})-\ell_{n-2}-2$ if $n\equiv0(\text{mod } 4)$  and $s_D(v_n)=2(\ell_0+\cdots+\ell_{n-2})+\ell_{n-1}$ if $n\equiv2(\text{mod } 4)$; and for any $2\le i\le n-1$,   $s_D(v_i)=2(\ell_0+\cdots+\ell_{i-2})+\ell_{i-1}+\ell_i$ if  $d_D^+(v_i)=0$, and   $s_D(v_i)=-2(\ell_0+\cdots+\ell_{i-3})-\ell_{i-2}-\ell_{i-1}-2$ if  $d_D^+(v_i)=2$.  When $n$ is odd, we see that $s_D(v_1)=\ell_0+\ell_1-1$;   $s_D(v_n)=-2(\ell_0+\cdots+\ell_{n-3})-\ell_{n-2}-2$ if $n\equiv3(\text{mod } 4)$  and $s_D(v_n)=2(\ell_0+\cdots+\ell_{n-2})+\ell_{n-1}$ if $n\equiv1(\text{mod } 4)$; and for any $2\le i\le n-1$,   $s_D(v_i)=2(\ell_0+\cdots+\ell_{i-2})+\ell_{i-1}+\ell_i$ if  $d_D^+(v_i)=0$,  and   $s_D(v_i)=-2(\ell_0+\cdots+\ell_{i-3})-\ell_{i-2}-\ell_{i-1}-2$ if  $d_D^+(v_i)=2$.  It can be easily checked that  for any $v_i, v_j\in V_R$ with $i\ne j$, $s_D(v_i)\neq s_D(v_j)$.
 \medskip

This completes the proof of Theorem~\ref{th2}.   \hfill\vrule height3pt width6pt depth2pt\\

 It would be nice if Theorem~\ref{th2} is true without assuming that $G$ is connected.  From the proof of Theorem~\ref{th2},  we obtain the following two results, where a component of a graph is \dfn{odd} if it has an odd number of vertices. \medskip

\begin{cor}
Let $G$ be a $2d$-regular graph, where $d\ge2$ is an integer. If $G$ has at most two odd components, then $G$ admits an antimagic orientation.
\end{cor}

\pf Let $G_1, G_2, \dots,  G_q$ be all the components of $G$, where $G_1$ is    the  odd component  (resp.  $G_1$ and $G_2$ are the odd components) of $G$  when $G$ has exactly one odd  component (resp. two odd components).  For each $i\in[q]$,   edges of $G_i$ are oriented as given in the proof of Theorem~\ref{th2},  and labeled  according to Algorithm~\ref{label2}.  Let $D$ be the resulting orientation of $G$. Clearly, the labeling of $D$ is   antimagic  if $G$ has at most one odd component. So we may assume that both $G_1$ and $G_2$ are odd. Let $v_1, \dots, v_n$ be the real vertices of an Euler tour of $G_1$,   and $u_1, u_2, \dots, u_m$ be the real vertices of an Euler tour of $G_2$.  From the proof of Theorem~\ref{th2}, no two vertices of $D$ has the same vertex-sum, except that $s_D(v_1)$ may be the same as $s_D(u_1)$. To avoid this, we relabel the edges on the paths $P_{1,2}$ and $P_{1,3}$   in the orientation of $G_1$ only as follows, where $P_{1,2}$,  $P_{1,3}$, $\ell_0$ and $\ell_1$ are defined as in the proof of Theorem~\ref{th2}:
assign the  numbers in $[\ell_1]$  to the edges of $P_{1, 3}$ in the  increasing order along  the orientation of  $P_{1,3}$, then assign the numbers in $\{\ell_1+1, \dots, \ell_1+\ell_0\}$ to  the edges of $P_{1,2}$ in the increasing order along the orientation of $P_{1,2}$.  One can easily check that  the resulting labeling of $D$ is antimagic.  \hfill\vrule height3pt width6pt depth2pt\medskip

\begin{cor}
Let $d\ge2$ be an integer. If every vertex of a connected graph $G$  has degree    $2d$ or $2d-2$,  then $G$ admits an antimagic orientation.
\end{cor}

\pf Let $G$ be a connected graph such that  every vertex of  $G$ has degree   $2d$ or $2d-2$. Let $C^*, C, D, D^*$ be  defined as in the proof of Theorem~\ref{th2}. Then $C^*$ contains  $d$ or $d-1$ copies of each vertex of $G$. From the proof of Theorem~\ref{th2}, we see that  $|s_D(u)-s_D(v)|\ge2$ for any two distinct vertices  $u, v$ in $D$ with $s_D(u)s_D(v)>0$; and  $s_{D^*}(u)=s_D(u)-(d-1) $ or $s_{D^*}(u)=s_D(u)-(d-2) $ for any $u$ in $D^*$. Thus $s_{D^*}(u)\ne s_{D^*}(v)$ for any two distinct vertices  $u, v$ in $D^*$.   Clearly, the labeling of $D^*$ is antimagic.  \hfill\vrule height3pt width6pt depth2pt\medskip

It seems hard to prove that if each of $G_1$ and $G_2$ has an antimagic orientation, then the disjoint union of $G_1$ and $G_2$  also has an antimagic orientation. But we know of no  counterexamples. With the support of Theorem~\ref{th1}, we believe the following is true.

\begin{conj}
Every  graph admits an antimagic orientation.
\end{conj}




\begin{thebibliography}{99}

\bibitem{Alon} N. Alon, G. Kaplan, A. Lev, Y. Roditty,  and R. Yuster, Dense graphs are antimagic, J Graph Theory 47 (2004) 297--309.
%
\vspace {-0.25cm}
%
\bibitem{CLPZ} F. Chang, Y-C. Liang, Z. Pan, and X. Zhu,  Antimagic labeling of
regular graphs,  J Graph Theory 82 (2016) 339--349.
%
\vspace {-0.25cm}
%
%
%
\bibitem{DC} D. Cranston, Regular bipartite graphs are antimagic, J Graph Theory 60 (2009) 173--182.
%
\vspace {-0.25cm}
%
\bibitem{DYZ} D. W. Cranston, Y-C. Liang, and X. Zhu, Regular graphs of odd degree are antimagic, J Graph Theory 80 (2015), 28--33.
%
\vspace {-0.25cm}
%
\bibitem{E} T.  Eccles,  {\small Graphs of large} linear size are antimagic,  J Graph Theory (2016) 81 236--261.
%
\vspace {-0.25cm}
%
\bibitem{JAG} J. A. Gallian, A dynamic survey of graph labeling, Electron J Combin DS6 (2016).
%
\vspace {-0.25cm}
%
\bibitem{NG} N. Hartsfield and G. Ringel, Pearls in Graph Theory, Academic Press, Boston, 1990, 108-109 (revised version, 1994).
%
\vspace {-0.25cm}
%
%
%
%
%
\bibitem{DTJ} D. Hefetz, T. M\"{u}tze, and J. Schwartz, On antimagic directed graphs, J Graph Theory 64 (2010) 219--232.
%
\vspace {-0.25cm}
%
%
%
%
%
\bibitem{SY} S. Shan and X. Yu, Antimagic orientation of biregular bipartite graphs, manuscript.  Available at   arXiv:1706.07336.

\end{thebibliography}
\end{document}